# APPROXIMATION AND SUPPORT THEOREM FOR A TWO SPACE-DIMENSIONAL WAVE EQUATION

ANNIE MILLET AND MARTA SANZ-SOLÉ

Abstract. We prove a characterization of the support of the law of the solution for a stochastic wave equation with two-dimensional space variable, driven by a noise white in time and correlated in space. The result is a consequence of an approximation theorem, in the convergence of probability, for equations obtained by smoothing the random noise. For some particular classes of coefficients, aproximation in the $L^p$ norm, for $p \geq 1$ is also proved.

## 1. Introduction and Preliminaries

In this paper we characterize the topological support of the law of the solution to the stochastic wave equation with two-dimensional spatial variable

$$(1.1) \qquad \begin{cases} \Big(\dfrac{\partial^2}{\partial t^2} - \Delta\Big)\, u(t,x) = \sigma\left(u(t,x)\right) F(dt,dx) + b\left(u(t,x)\right) \\ u(0,x) = u_0(x) \\ \dfrac{\partial u}{\partial t}\ (0,x) = v_0(x) \end{cases}$$

$(t,x) \in [0,\infty[\times\mathbb{R}^2$.

Here $F(t,x)$ is a generalized Gaussian field with covariance

$$(1.2) \qquad E\Big(F(t,x)\ F(s,y)\Big) = \delta(t-s)\ f(\,|x-y|\,),$$

where $\delta$ denotes the Dirac delta function and $f$ is a continuous function on $]0,\infty[$ satisfying an integrability condition made precise later (see (C1)). In addition the functional $J : \mathcal{D}(\mathbb{R}^3) \times \mathcal{D}(\mathbb{R}^3) \longrightarrow \mathbb{R}$ given by

$$J(\varphi,\psi) = \int_0^\infty dt \int_{\mathbb{R}^2} dx \int_{\mathbb{R}^2} dy\ \varphi(t,x)\ f(|x-y|)\ \Psi(t,y)$$

is positive definite. With these hypotheses the process $\{F(t,x), (t,x) \in [0,\infty[\times\mathbb{R}^2\}$ exists.

We consider the mild formulation of (1.1). That means, let

$$S(t,x) = \frac{1}{2\pi}\ (t^2 - |x|^2)^{-\frac{1}{2}}\ 1_{\{|x|<t\}}\,,$$

Part of this work has been done during a visit of the authors to the Mathematical Sciences Research Institute at Berkeley. We would like to thank this Center for its hospitality and support. Both authors have also been partially supported by the grant ERBF MRX CT960075A from the European Union. The second named author has been supported by the grant PB96-0088 from the *Subdirección General de Formación y Promoción del Conocimiento.*

then a solution of (1.1) is a stochastic process $\{u(t,x),\ (t,x)\in[0,\infty[\times\mathbb{R}^2\}$ satisfying

$$
\begin{aligned}
u(t,x) &= \int_{\mathbb{R}^2} S(t,x-y)\ v_0(y)\,dy + \frac{\partial}{\partial t}\Big(\int_{\mathbb{R}^2} S(t,x-y)\ u_0(y)\,dy\Big)\\
&\quad + \int_0^t \int_{\mathbb{R}^2} S(t-s,\ x-y)\ [\sigma\,(u(s,y))\ F(ds,dy) + b\,(u(s,y))\,ds\,dy].
\end{aligned}
\tag{1.3}
$$

Consider the following set of assumptions on the elements defining (1.3):

(C1) There exists $\beta\in(0,1)$, $r_0>0$, such that for $0<t<r_0$

$$
\int_0^t r\,f(r)\ \ln\Big(1+\frac{t}{r}\Big)dr \le C\,t^\beta.
$$

(C2) $u_0:\mathbb{R}^2\longrightarrow\mathbb{R}$ is of class $\mathcal{C}^1$ and bounded with $\frac{\beta}{2(1+\beta)}$-Hölder continuous partial derivatives, $v_0:\mathbb{R}^2\longrightarrow\mathbb{R}$ is such that $|v_0|+|\nabla u_0|\in L^{q_0}$ for some $q_0\in]4,\infty]$.

(C3) $\sigma,b:\mathbb{R}\longrightarrow\mathbb{R}$ are globally Lipschitz functions.

Fix $T>0$ and a compact set $K\subset\mathbb{R}^2$. In the previous work Millet and Sanz-Solé (1997), we have proved that, assuming (C1), (C2) and (C3), equation (1.3) has a unique solution. Moreover, the trajectories of $u$ are $\gamma$-Hölder continuous in $(t,x)\in[0,T]\times K$ for every $\gamma\in\Big(0,\ \frac{\beta}{2(1+\beta)}\Big)$.

Let $H$ denote the completion of the inner-product space consisting of measurable functions $h:[0,T]\times\mathbb{R}^2\longrightarrow\mathbb{R}$ such that

$$
\int_0^T ds\int_{\mathbb{R}^2} dx\int_{\mathbb{R}^2} dy\ |h|\,(s,x)\ f(\,|x-y|\,)\ |h|\ (s,y) < +\infty\,,
$$

endowed with the inner product

$$
\langle h_1,h_2\rangle_H := \int_0^T ds\int_{\mathbb{R}^2} dx\int_{\mathbb{R}^2} dy\ h_1(s,x)\ f(\,|x-y|\,)\ h_2(s,y)\,.
$$

For any $h\in H$, let $\{\Phi^h(t,x),\ (t,x)\in[0,\infty)\times\mathbb{R}^2\}$ be the solution of

$$
\begin{aligned}
\Phi^h(t,x) &= \int_{\mathbb{R}^2} S(t,x-y)\,v_0(y)dy + \frac{\partial}{\partial t}\left(\int_{\mathbb{R}^2} S(t,x-y)\,u_0(y)\,dy\right)\\
&\quad + \int_0^t ds\int_{\mathbb{R}^2} dy\int_{\mathbb{R}^2} dz\ S(t-s,\ x-y)\ \sigma(\Phi^h(s,y))\,f(\,|y-z|\,)\,h(s,z)\\
&\quad + \int_0^t\int_{\mathbb{R}^2} S(t-s,\ x-y)\ b(\Phi^h(s,y))\,ds\,dy\,.
\end{aligned}
\tag{1.4}
$$

Set

$$
\|\varphi\|_{\gamma,K} = \sup_{\substack{t\in[0,T]\\ x\in K}} |\varphi(t,x)| + \sup_{\substack{t,t'\in[0,T]\\ x,x'\in K\\ t\neq t',\ x\neq x'}} \frac{|\varphi(t,x)-\varphi(t',x')\,|}{(\,|t-t'|+|x-x'|\,)^\gamma}\,.
\tag{1.5}
$$

We prove in Section 2 that the support of the law of $\{u(t,x),\ (t,x)\in[0,T]\times K\}$ is the closure with respect to the norm $\|\cdot\|_{\gamma,k}$ of the set of functions $\{\Phi^h,\ h\in H\}$, where $\{\Phi^h(t,x),\ (t,x)\in[0,T]\times K\}$ is the solution of (1.4). The proof is based on an approximation result for equations more general than (1.3) constructed by smoothing the random noise $F(t,x)$. We refer the reader to Millet and Sanz-Solé (1994a), Millet and Sanz-Solé (1994b) and Bally, Millet and Sanz-Solé (1995) for a presentation of the method an applications to stochastic differential and stochastic partial differential equations.

In the framework of stochastic partial differential equations, the regularization of the noise rises up technical difficulties connected with the explosion of the corresponding integral (see for instance Bally, Millet and Sanz-Solé (1995)). This problem does not appear here because of the following reasons. The noise $F$ is smoother than space-time white noise. On the other hand the integrability condition (C1) and Lemma A.1 in Millet and Sanz-Solé (1997) yield

$$\mu(t) \leq C\ t^{\beta+1}$$

(see (4.2) and (4.11) ). This fact prevents explosions, as in made explicit in the proofs. We now introduce some preliminaries and notations used along the paper.

Let $\tilde{H}$ be the completion of the inner-product space of measurable functions $\varphi : \mathbb{R}^2 \longrightarrow \mathbb{R}$ such that $\int_{\mathbb{R}^2} dx \int_{\mathbb{R}^2} dy\ |\varphi(x)|\ f(|x-y|)\ |\varphi(y)| < +\infty$ endowed with the inner product

$$\langle \varphi, \psi \rangle_{\tilde{H}} := \int_{\mathbb{R}^2} dx \int_{\mathbb{R}^2} dy\ \varphi(x)\ f(\,|x-y|\,)\ \psi(y)\,.$$

Let $\{e_j,\ j \in \mathbb{N}\}$ be a complete orthonormal system of $\tilde{H}$ that is supposed to be fixed in the sequel. Define

$$W_j(t) = \int_0^t \int_{\mathbb{R}^2} e_j(x)\ F(ds, dx),\ j \in \mathbb{N},\ t \in [0,T]\,. \tag{1.6}$$

Clearly $\{W_j,\ j \in \mathbb{N}\}$ is a sequence of independent Brownian motions.

Let $\mathcal{H}$ be the separable Hibert space of functions $k : [0,T] \longrightarrow \mathbb{R}^{\mathbb{N}}$ such that $\int_0^T \sum_{j=1}^{\infty} |k_j(s)|^2\, ds < \infty$ endowed with the inner product

$$\langle k, \bar{k} \rangle_{\mathcal{H}} = \int_0^T \sum_{j=1}^{\infty} k_j(s)\ \bar{k}_j(s)\ ds\,.$$

Notice that the mapping

$$\begin{aligned} \mathcal{J}:\ H &\longrightarrow \mathcal{H} \\ \varphi &\longmapsto \Big( \int_0^{\cdot} \langle \varphi(s,*),\, e_j(*) \rangle_{\tilde{H}}\, ds \Big)_{j \in \mathbb{N}} \end{aligned} \tag{1.7}$$

provides an isometry between $H$ and $\mathcal{H}$.

Let $\mathcal{F}_t = \sigma\Big( F(\,[0,s] \times A);\ 0 \leq s \leq t\,,\ A \in \mathcal{B}(\mathbb{R}^2) \Big),\ t \geq 0$. For any $\mathcal{F}_t$-predictable process $\varphi \in L^2(\Omega; H)$ we have

$$\int_0^t \int_{\mathbb{R}^2} \varphi(s,x)\ F(ds,dx) = \sum_{j=1}^{\infty} \int_0^t \langle \varphi(s,*),\, e_j(*) \rangle_{\tilde{H}}\ W_j(ds)\,, \tag{1.8}$$

$t \in [0,T]$, so that the stochastic integral with respect to the martingale measure $F$ can also be viewed as a stochastic integral with respect to the Gaussian process $\{W(k),\ k \in \mathcal{H}\}$ defined by

$$W(k) = \sum_{j=1}^{\infty} \int_0^T k_j(s)\, W_j(ds)\,. \tag{1.9}$$

We introduce smooth approximations of $F$ constructed as follows. Fix $n \in \mathbb{N}$ and let $\mathcal{P}_n$ be the partition of $[0,T]$ determined by $\frac{iT}{2^n},\ i = 0,1,\ldots,2^n$. Denote by $\Delta_i$ the interval $\Big[\frac{iT}{2^n}, \frac{(i+1)T}{2^n}\Big)$ and by $|\Delta_i|$ its length. We write $W_j(\Delta_i)$ for the increment $W_j\Big(\frac{(i+1)T}{2^n}\Big) - W_j\big(\frac{iT}{2^n}\big)$, $i = 0,\ldots,2^n-1$. Define

$$W^n = \Big( W_j^n = \int_0^{\cdot} \dot{W}_j^n(s)ds\,,\ j \in \mathbb{N}\Big) \tag{1.10}$$

where $W_j^n = 0$ if $j > n$, and for $1 \le j \le n$,

$$\dot{W}_j^n(t) = \begin{cases} \displaystyle\sum_{i=1}^{2^n} 2^n\, T^{-1}\, W_j(\Delta_{i-1})\ 1_{\Delta_i}(t) & \text{if} \quad t \in [2^{-n}T, T]\,, \\ 0 & \text{if} \quad t \in [0,\, 2^{-n}T)\,. \end{cases} \tag{1.11}$$

Set

$$\omega^n(t,x) = \sum_{j\in\mathbb{N}} \dot{W}_j^n(t)\ e_j(x)\,. \tag{1.12}$$

It is easy to check that a.s., $\omega^n \in H$ and, more precisely,

$$\|\omega^n\|_{L^p(\Omega;H)} \le C\, n^{\frac12}\ 2^{\frac{n}{2}},\ \forall p \in [1,\,\infty)\,. \tag{1.13}$$

Indeed, fix $p \in [2,\infty)$; then

$$\begin{aligned} \|\omega^n\|^p_{L^p(\Omega;H)} &= E\Big( \Big| \sum_{j=1}^{n} \sum_{i=1}^{2^n} 2^{2n}\ T^{-2}\ |\Delta_i|\ W_j(\Delta_{i-1})^2\ \Big|^{\frac{p}{2}} \Big) \\ &\le C\ n^{\frac{p}{2}}\ 2^{\frac{np}{2}}\,. \end{aligned}$$

Moreover, for any $0 \le t_1 \le t_2 \le T$, similar computations imply

$$\|\omega^n\ 1_{[t_1,t_2]}\|_{L^p(\Omega;H)} \le C\ n^{\frac12}\ 2^{\frac{n}{2}}\ |t_2 - t_1|^{\frac12}. \tag{1.14}$$

Let $(\bar\Omega, \bar{\mathcal{F}}, \bar P)$ be the canonical space associated with a standard Brownian motion. Denote by $(\Omega, \mathcal{F}, P)$ the product space $(\bar\Omega^{\mathbb{N}}, \bar{\mathcal{F}}^{\otimes\mathbb{N}}, \bar P^{\otimes\mathbb{N}})$, which will be our reference probability space.
Set $\bar k(t) = \int_0^t k(s)\, ds$ for $k \in \mathcal{H}$. For any integer $n \ge 1$, let $T_n^k$ denote the transformation of $\Omega$ defined by

$$T_n^k(w) = w + \bar k - w^n\,. \tag{1.15}$$

Notice that $T_n^k(w) = w + \int_0^{\cdot} \varphi_n(s,w)\, ds$, where $\{\varphi_n(t,w),\ t\in[0,T]\}$ is an $\mathcal{H}$-valued process adapted to the filtration generated by $\{W_j(t),\ t\in[0,1],\ j\in\mathbb{N}\}$. Therefore, by Girsanov's theorem, $P_\circ(T_n^k)^{-1} \ll P$. This fact will be used in the proof of Theorem 2.1.
The paper is organized as follows. In Section 2 we prove the characterization of the support by means of an approximation in probability. In Section 3 we prove approximations in $L^p$-norm under stronger hypotheses on the coefficients. As usual, all constants are denoted by $C$, independently of their values.

## 2. Approximation in probability and support theorem

The purpose of this section is to prove the following result.

**Theorem 2.1.** *Assume (C1) to (C3), fix a compact set $K\subset\mathbb{R}^2$ and let $\{u(t,x),\ t\in[0,T],\ x\in K\}$ be the solution of (1.3). Then for any $\gamma\in\Big(0,\ \frac{\beta}{2(1+\beta)}\Big)$ the topological support of the law of $u$ in the space $\mathcal{C}^\gamma([0,T]\times K)$ of $\gamma$-Hölder continuous functions in $(t,x)$ is given by the closure in $\mathcal{C}^\gamma([0,T]\times K)$ of the set of functions $\{\Phi^h,\ h\in H\}$, where $\{\Phi^h(t,x),\ t\in[0,T],\ x\in K\}$ is the solution of (1.4).*

The proof of Theorem 2.1 is a consequence of an approximation result, in the convergence in probability, for an equation more general that (1.3).
More precisely, let us introduce the following hypothesis

(C3') The coefficients $A,B,D,b:\mathbb{R}\longrightarrow\mathbb{R}$ are globally Lipschitz functions.

Then we consider the evolution equations

$$\begin{aligned} X_n(t,x) &= X^0(t,x) + \int_0^t\int_{\mathbb{R}^2} S(t-s,\,x-y)\, A\left(X_n(s,y)\right) F(ds,\,dy) \\ &\quad + \langle S\left(t-\cdot,\,x-*\right) B(X_n(\cdot,*)),\ \omega^n\rangle_H + \ \langle S\left(t-\cdot,\,x-*\right) D(X_n(\cdot,*)),\ h\rangle_H \\ &\quad + \int_0^t\int_{\mathbb{R}^2} S(t-s,\,x-y)\, b\left(X_n(s,y)\right) ds\, dy\,, \end{aligned} \tag{2.1}$$

$$\begin{aligned} X(t,x) &= X^0(t,x) + \int_0^t\int_{\mathbb{R}^2} S(t-s,\,x-y)\ (A+B)\ \left(X(s,y)\right) F(ds,\,dy) \\ &\quad + \langle S\left(t-\cdot,\,x-*\right) D(X(\cdot,*)),\ h\rangle_H \\ &\quad + \int_0^t\int_{\mathbb{R}^2} S(t-s,\,x-y)\, b\left(X(s,y)\right) ds\, dy\,, \end{aligned} \tag{2.2}$$

where $n\ge 1,\ A,B,D,b$ satisfy (C3'), $h\in H$, $\omega^n$ is defined in (1.12) and

$$X^0(t,x) = \int_{\mathbb{R}^2} S(t,\,x-y)\ v_0(y)\, dy + \frac{\partial}{\partial t}\Big(\int_{\mathbb{R}^2} S(t,\,x-y)\ u_0(y)\, dy\Big). \tag{2.3}$$

Our aim is to prove the following:

**Proposition 2.2.** *Assume (C1), (C2) and (C3'). For any* $\gamma \in \Big(0, \frac{\beta}{2(1+\beta)}\Big)$, $\eta > 0$ *and every compact set* $K \subset \mathbb{R}^2$,

$$\lim_{n\to\infty} P\left( \|X_n - X\|_{\gamma,K} > \eta\right) = 0\,, \tag{2.4}$$

*where* $\|\cdot\|_{\gamma,K}$ *has been defined in* (1.5).

We at first show that Theorem 2.1 is an easy consequence of this Proposition.

*Proof of Theorem 2.1.* Assume that Proposition 2.2 has been proved. For $n \geq 1$, set

$$\begin{aligned} u_n(t,x) &= X^0(t,x) + \langle S(t-\cdot,\, x-*)\; \sigma(u_n(\cdot,*)),\; \omega^n\rangle_H \\ &\quad + \int_0^t \int_{\mathbb{R}^2} S\,(t-s,\, x-y)\; b(u_n(s,y))\, ds\, dy\,, \end{aligned} \tag{2.5}$$

$$\begin{aligned} v_n(t,x) &= X^0(t,x) + \int_0^t \int_{\mathbb{R}^2} S(t-s,\, x-y)\; \sigma(v_n(s,y))\; F(ds,\, dy) \\ &\quad + \langle S(t-\cdot,\, x-*)\; \sigma(v_n(\cdot,*)),\; h-\omega^n\rangle_H \\ &\quad + \int_0^t \int_{\mathbb{R}^2} S(t-s,\, x-y)\; b(v_n(s,y))\, ds\, dy\,. \end{aligned} \tag{2.6}$$

Clearly, equations (2.5) and (2.6) are particular cases of (2.1) while equations (1.3) and(1.4) are particular cases of (2.2), obtained by choosing $A = D = 0,\ B = \sigma$ and $A = D = \sigma,\ B = -\sigma$, respectively.
Moreover, $u_n = \Phi^{\omega^n}$. Given $h \in H$, set $k = \mathcal{J}(h)$, where $\mathcal{J}$ is the isometry defined in (1.7). Then, by (1.8), equation (2.6) can be rewritten as follows

$$\begin{aligned} v_n(t,x) &= X^0(t,x) + \sum_{j=1}^{\infty} \int_0^t \langle S(t-s,\, x-*)\; \sigma(v_n(s,*)),\; e_j\rangle_{\tilde H}\; W_j(ds) \\ &\quad + \sum_{j=1}^{\infty} \int_0^t \langle S(t-s,\, x-*)\; \sigma(v_n(\cdot,*)),\; e_j\rangle_{\tilde H}\, (k_j(s) - \dot W_j^n(s))\, ds \\ &\quad + \int_0^t \int_{\mathbb{R}^2} S(t-s,\, x-y)\; b(\sigma_n(s,y))\; ds\, dy\,, \end{aligned}$$

with $\dot W_j^n$ defined in (1.11). Therefore, $v^n = u \circ T_n^k$, where $T_n^k$ is the absolutely continuous transformation on $\Omega$ defined by (1.15).
The convergence (2.4) implies for any $\eta > 0$,

$$\begin{aligned} \lim_{n\to\infty} P(\, \|\Phi^{\omega^n} - u\|_{\gamma,K} > \eta) &= 0, \\ \lim_{n\to\infty} P(\, \|u(T_n^h) - \Phi^h\|_{\gamma,K} > \eta) &= 0\,. \end{aligned}$$

These two convergences yield the characterization of the support stated in Theorem 2.1 (see, for instance Bally, Millet and Sanz-Solé (1995)). □

Before proving Proposition 2.2, we should address the question of existence and uniqueness of solution for the equations (2.1) and (2.2), respectively.

As in Millet and Sanz-Solé (1997), the classical Picard iteration scheme provides the existence of a unique solution of (2.2) in the space $\mathcal{C}^\gamma([0,T]\times K),\ \gamma\in\Big(0,\frac{\beta}{2(1+\beta)}\Big)$. Moreover,

$$\sup_{0\le t\le T}\ \sup_{x\in\mathbb{R}^2} E(\,|X(t,x)|^p) < \infty\,,\quad p\in[1,\infty)\,. \tag{2.7}$$

This method does not seem appropiate for equation (2.1), due to the term involving $\omega^n$ which has an unbounded $H$-norm. For this reason, we first localize $\omega^n$ as follows. For any positive integer $n,\ M\in\mathbb{R}_+$ and $t\in[0,T]$, set

$$A_{n,M}(t) = \Big\{ \sup_{1\le j\le n}\ \sup_{0\le i\le 2^{-n}([2^n tT^{-1}]-1)^+} 2^n\,|W_j(\Delta_i)\,| \le M\Big\} \tag{2.8}$$

and

$$\omega^{n,M}(t,x) = \omega^n(t,x)\ 1_{A_{n,M}(t)}.$$

Notice that

$$\sup_{0\le t\le T} \|\omega^{n,M}(t,*)\|_{\tilde H} \le M\ \sqrt{n}\,. \tag{2.9}$$

Consider the evolution equation

$$\begin{aligned} X_{n,M}(t,x) \;&=\; X^0(t,x) + \int_0^t\int_{\mathbb{R}^2} S(t-s,\ x-y)\ A(X_{n,M}(s,y))\ F(ds,dy)\\ &\quad + \langle S(t-\cdot,\ x-*)\ B(X_{n,M}(\cdot,*)),\ \omega^{n,M}\rangle_H\\ &\quad + \langle S(t-\cdot,\ x-*)\ D(X_{n,M}(\cdot,*)),\ h\,\rangle_H\\ &\quad + \int_0^t\int_{\mathbb{R}^2} S(t-s,\ x-y)\ b(X_{n,M}(s,y))\ ds\,dy\,. \end{aligned}$$

As in Millet and Sanz-Solé (1997), Picard's iteration scheme provides the existence and uniqueness of the solution to this equation. For any $(t,x)\in[0,1]\times\mathbb{R}^2$, define

$$X_n(t,x) = X_{n,M}(t,x)\quad \text{on}\quad A_{n,M}(T)\,.$$

For fixed $n$ and $M$, the sets $(A_{n,M}(t))_{t\in[0,T]}$ are decreasing in $t$. Therefore, a standard argument based on the local property of stochastic integrals implies that this definition is consistent and, since $P(\cup_{M\ge1} A_{n,M}(1)) = 1$ for every integer $n$, this provides the existence and uniqueness of solution to equation (2.1).
The proof of Proposition 2.2 relies on a localization procedure and on Lemma 4.1. We start by giving the ingredients which are needed in the localization.
Fix $\alpha > (2\ \ln\ 2)^{\frac12}$ and for every $n>0$, set

$$M(n) = \alpha\ 2^{\frac n2}\ n^{\frac12}\,, \tag{2.10}$$

and

$$A_n(t) = A_{n,M(n)}\ (t)\,. \tag{2.11}$$

**Lemma 2.3.** *The following convergence holds:*

$$\lim_{n\to\infty} P\Big(A_n(T)^c\Big) = 0\,.$$

*Proof*: Let $Z$ denote a $N(0,1)$ random variable. Then

$$\begin{aligned} P\Big(A_n(T)^c\Big) &\leq n\, 2^n\, P\Big(|Z| > 2^{-\frac{n}{2}}\, M(n)\Big) \\ &\leq C\, n\, 2^n\, \frac{2^{\frac{n}{2}}}{M(n)}\, \exp\left(-\frac{2^{-n}\, M(n)^2}{2}\right) \\ &= C\, \sqrt{n}\, \exp\left(-n\Big(\frac{\alpha^2}{2} - \ln 2\Big)\right) \xrightarrow[n\to\infty]{} 0\,. \end{aligned}$$

□

**Remark 2.4.** Due to (2.9), on the set $A_n(T)$, we have:

$$\|\omega^n\|_H \leq C\; n\; 2^{\frac{n}{2}} \tag{2.12}$$

and, for any $0 \leq t \leq t' \leq T$, on $A_n(t')$ we have:

$$\|\omega^n\; 1_{[t,t']}\|_H \leq C\; n\; 2^{\frac{n}{2}}\, |t'-t|^{\frac{1}{2}}\,. \tag{2.13}$$

In particular, if $[t,t'] \subset \Delta_i$ for some $i = 0, \ldots, 2^n - 1$, on $A_n(t')$ it holds:

$$\|\omega^n\; 1_{[t,t']}\|_H \leq C\; n\,. \tag{2.14}$$

Our next purpose is to check that the sequence of processes $Y_n(t,x) := X_n(t,x) - X(t,x)$ , $n \geq 1$ satisfies the requirements of Lemma 4.1.
To this end, we introduce some notations and prove several Lemmas. For any $n \geq 1$, $t \in [0,T]$, set

$$\begin{aligned} \underline{t}_n &= \max\{k\, 2^{-n}\, T;\; k = 1, \ldots, 2^n - 1 : k\, 2^{-n}\, T \leq t\}\,, \\ t_n &= \big(\underline{t}_n - 2^{-n}\, T\big) \vee 0\,, \end{aligned} \tag{2.15}$$

$$\begin{aligned} X_n^-(t,x) &= X^0(t,x) + \int_0^{t_n}\int_{\mathbb{R}^2} S(t-s,\, x-y)\; A(X_n(s,y))\, F(ds,\, dy) \\ &\quad + \langle S(t-\cdot,\, x-*)\; B(X_n(\cdot,*))\; 1_{[0,t_n]}(\cdot)\,,\; \omega^n\rangle_H \\ &\quad + \langle S(t-\cdot,\, x-*)\; D(X_n(\cdot,*))\; 1_{[0,t_n]}(\cdot)\,,\; h\rangle_H \\ &\quad + \int_0^{t_n}\int_{\mathbb{R}^2} S(t-s,\, x-y)\; b(X_n(s,y))\, ds\, dy\,, \end{aligned} \tag{2.16}$$

$$\begin{aligned} X^-(t,x) &= X^0(t,x) + \int_0^{t_n}\int_{\mathbb{R}^2} S(t-s,\, x-y)\; (A+B)\,(X(s,y))\, F(ds,\, dy) \\ &\quad + \langle S(t-\cdot,\, x-*)\; D(X(\cdot,*))\; 1_{[0,t_n]}(\cdot)\,,\; h\rangle_H \\ &\quad + \int_0^{t_n}\int_{\mathbb{R}^2} S(t-s,\, x-y)\; b(X(s,y))\, ds\, dy\,. \end{aligned} \tag{2.17}$$

To lighten the notations, we do not write explicitely the fact that the process $X^-$ depends on $n$. In the sequel $\|\;\|_p$ denotes the $L^p(\Omega)$-norm.

**Lemma 2.5.** *Suppose that conditions (C1), (C2) and (C3') hold. Then, for any $p \in [1,\infty)$ and every integer $n \geq 1$,*

$$\sup_{(s,x)\in[0,T]\times\mathbb{R}^2} \|X(s,x) - X^-(s,x)\|_p \leq C\; 2^{-n\,\frac{\beta+1}{2}}.$$

*Proof*: Set $\|X(t,x) - X^-(t,x)\|_p^p \leq C \sum_{i=1}^3 V_i(t,x)$, with

$$\begin{aligned} V_1(t,x) &= E\left(\Big|\int_{t_n}^t \int_{\mathbb{R}^2} S(t-s,\, x-y)\; (A+B)\; (X(s,y))\; F(ds,\, dy)\Big|^p\right), \\ V_2(t,x) &= E\Big(|\langle S(t-\cdot,\, x-*)\, D(X(\cdot,*))\; 1_{(t_n,t]}(\cdot),\; h\rangle_H|^p\Big), \\ V_3(t,x) &= E\left(\Big|\int_{t_n}^t \int_{\mathbb{R}^2} S(t-s,\, x-y)\; b\,(X(s,y))\, ds\, dy\Big|^p\right). \end{aligned}$$

Burkholder's and Hölder's inequalities, (2.7) and (4.11) yield

$$V_1(t,x) \leq C\; \mu(t-t_n)^{\frac{p}{2}} \left(1 + \sup_{(t,x)\in[0,T]\times\mathbb{R}^2} E\Big(|X(t,x)|^p\Big)\right) \leq C\; 2^{-n(\beta+1)\frac{p}{2}},$$

with $\mu(t-t_n)$ given by (4.2).
Schwarz's and Hölder's inequalities imply

$$V_2(t,x) \leq C\; \|h\|_H^p\; \mu(t-t_n)^{\frac{p}{2}} \left(1 + \sup_{(t,x)\in[0,T]\times\mathbb{R}^2} E\Big(|X(t,x)|^p\Big)\right) \leq C\; 2^{-n(\beta+1)\,\frac{p}{2}}.$$

Finally, Hölder's inequality implies for $\nu(t)$ defined by (4.3):

$$\begin{aligned} V_3(t,x) &\leq C\; \Big(\int_{t_n}^t \int_{\mathbb{R}^2} S(t-s,\, x-y)\, ds\, dy\Big)^p \left(1 + \sup_{\substack{0\leq t\leq T\\ x\in\mathbb{R}^2}} E\Big(|X(t,x)|^p\Big)\right) \\ &\leq C\; \nu(t-t_n)^p \leq C\; 2^{-2np}, \end{aligned}$$

which completes the proof of the Lemma. □

Consider the Picard iteration scheme associated with (2.1):

$$\begin{aligned} X_n^0(t,x) &= X^0(t,x) \quad \text{and for } k \geq 0\,, \\ X_n^{k+1}(t,x) &= X^0(t,x) + \int_0^t \int_{\mathbb{R}^2} S(t-s,\, x-y)\; A(X_n^k(s,y))\; F(ds,\, dy) \\ &\quad + \langle S(t-\cdot, x-*)\; B(X_n^k(\cdot,*)),\; \omega^n\rangle_H + \langle S(t-\cdot, x-*)\; D(X_n^k(\cdot,*)),\; h\rangle_H \\ &\quad + \int_0^t \int_{\mathbb{R}^2} S(t-s,\, x-y)\; b(X_n^k(s,y))\; ds,\, dy\,. \end{aligned} \tag{2.18}$$

For any $0 \le r \le t \le T$ and every integer $k \ge 0$ set

$$\begin{aligned} X_n^0(t,r;x) &= X^0(t,x)\,, \\ X_n^{k+1}(t,r;x) &= X^0(t,x) + \int_0^r \int_{\mathbb{R}^2} S(t-s,\, x-y)\; A(X_n^k(s,y))\; F(ds,\, dy) \\ &\quad + \langle S(t-\cdot, x-*)\; B(X_n^k(\cdot,*))\; 1_{[0,r]}\,(\cdot)\,,\; \omega^n\rangle_H \\ &\quad + \langle S(t-\cdot, x-*)\; D(X_n^k(\cdot,*))\; 1_{[0,r]}\,(\cdot)\,,\; h\rangle_H \\ &\quad + \int_0^r \int_{\mathbb{R}^2} S(t-s,\, x-y)\; b(X_n^k(s,y))\; ds\, dy\,, \qquad (2.19) \\ \bar{X}_n^{k+1}(t,x) &= X_n^{k+1}\,(t, t_n; x)\,. \qquad (2.20) \end{aligned}$$

Notice that $X_n^k(t,t;x) = X_n^k(t,x)$.

**Lemma 2.6.** *Assume (C1), (C2) and (C3'). Then, for every $p \in [1,\infty)$, $t \in [0,T]$, $k \ge 1$, $n \ge 1$,*

$$\begin{aligned} \sup_{(s,y)\in[0,t]\times\mathbb{R}^2} E\left(|X_n^k(s,y) - \bar{X}_n^k(s,y)|^p\; 1_{A_n(s)}\right) &\le C\; n^p\; 2^{-n(1+\beta)\,\frac{p}{2}} \\ &\quad \times \Big[1 + \sup_{(s,y)\in[0,t]\times\mathbb{R}^2} E\big(\,|X_n^{k-1}(s,y)|^p\; 1_{A_n(s)}\big)\Big]\,. \qquad (2.21) \end{aligned}$$

*and*

$$\begin{aligned} \sup_{(s,y)\in[0,t]\times\mathbb{R}^2} E\left(|X_n(s,y) - X_n^-(s,y)|^p\; 1_{A_n(s)}\right) &\le C\; n^p\; 2^{-n(1+\beta)\,\frac{p}{2}} \\ &\quad \times \Big[1 + \sup_{(s,y)\in[0,t]\times\mathbb{R}^2} E\big(\,|X_n(s,y)|^p\; 1_{A_n(s)}\big)\Big]\,. \qquad (2.22) \end{aligned}$$

*Proof*. Consider the decomposition

$$E\Big(|X_n^k(t,x) - \bar{X}_n^k(t,x)|^p\; 1_{A_n(t)}\Big) \le C\; \sum_{i=1}^4\; T_n^{k,i}(t,x)\,, \qquad (2.23)$$

with

$$\begin{aligned} T_n^{k,1}(t,x) &= E\left(\Big|\int_{t_n}^t \int_{\mathbb{R}^2} S(t-s,\, x-y)\; A(X_n^{k-1}(s,y))\; F(ds,\, dy)\Big|^p\; 1_{A_n(t)}\right), \\ T_n^{k,2}(t,x) &= E\left(\Big|\langle S(t-\cdot,\, x-*)\; B(X_n^{k-1}(\cdot,\, *))\; 1_{(t_n,t]}\,(\cdot),\; \omega^n\rangle_H\Big|^p\, 1_{A_n(t)}\right), \\ T_n^{k,3}(t,x) &= E\left(\Big|\langle S(t-\cdot,\, x-*)\; D(X_n^{k-1}(\cdot,\, *))\; 1_{(t_n,t]}\,(\cdot),\; h\rangle_H\Big|^p\, 1_{A_n(t)}\right), \\ T_n^{k,4}(t,x) &= E\left(\Big|\int_{t_n}^t \int_{\mathbb{R}^2} S(t-s,\, x-y)\; b(X_n^{k-1}(s,y))\; ds\, dy\Big|^p\; 1_{A_n(t)}\right). \end{aligned}$$

The local property of stochastic integrals, the inclusion $A_n(s) \supset A_n(t)$ for $s \le t$, Burkholder's and Hölder's inequalities and (4.11) yield

$$\begin{aligned} T_n^{k,1}(t,x) &\le C\,\mu\,(t-t_n)^{\frac{p}{2}}\,\Big[1+\sup_{(s,y)\in[0,t]\times\mathbb{R}^2} E(\,|X_n^{k-1}(s,y)\,|^p\,1_{A_n(s)})\Big] \\ &\le C\,2^{-n\,(\beta+1)\,\frac{p}{2}}\,\Big[1+\sup_{(s,y)\in[0,t]\times\mathbb{R}^2} E(\,|X_n^{k-1}(s,y)\,|^p\,1_{A_n(s)})\Big]. \end{aligned} \tag{2.24}$$

Schwarz's and Hölder's inequalities, (2.14) and (4.11) imply

$$\begin{aligned} T_n^{k,2}(t,x) &\le E\Big[\|\omega^n\,1_{(t_n,t]}\,1_{A_n(t)}\|_H^p \\ &\qquad\times\|S(t-\cdot,\,x-*)\,B(X_n^{k-1}\,(\cdot,*))\,1_{(t_n,t]}\,(\cdot)\,1_{A_n(t)}\|_H^p\Big] \\ &\le C\,n^p\,2^{-n(1+\beta)\,\frac{p}{2}}\,\Big[1+\sup_{(s,y)\in[0,t]\times\mathbb{R}^2} E(\,|X_n^{k-1}(s,y)\,|^p\,1_{A_n(s)})\Big]. \end{aligned} \tag{2.25}$$

Similarly, using (4.3) for the last inequality, we have

$$T_n^{k,3}(t,x) \le C\,\|h\|_H^p\,2^{-n(1+\beta)\,\frac{p}{2}}\,\Big[1+\sup_{(s,y)\in[0,t]\times\mathbb{R}^2} E(\,|X_n^{k-1}(s,y)|^p\,1_{A_n(s)})\Big], \tag{2.26}$$

$$T_n^{k,4}(t,x) \le C\,2^{-2np}\,\Big[1+\sup_{(s,y)\in[0,t]\times\mathbb{R}^2} E(\,|X_n^{k-1}(s,y)|^p\,1_{A_n(s)})\Big]. \tag{2.27}$$

Thus (2.23)–(2.27) conclude the proof of (2.21).
Using the arguments in the proof of Theorem 1.2 Millet and Sanz-Solé (1997), we obtain

$$\lim_{k\to\infty}\sup_{(s,x)\in[0,t]\times\mathbb{R}^2} E\Big[\Big(|X_n^k(s,x)-X_n(s,x)|^p+|\bar{X}_n^k(s,x)-X_n^-(s,x)|^p\Big)\,1_{A_n(s)}\Big]=0\,. \tag{2.28}$$

Therefore, (2.21) and (2.28) yield (2.22). □

We now prove the convergence of $X_n^-(s,y)$ to $X_n(s,y)$.

**Lemma 2.7.** *Assume (C1), (C2) and (C3'). Then for any* $p\in[1,+\infty)$,

$$\sup_{n\ge1}\sup_{(t,x)\in[0,T]\times\mathbb{R}^2} E\Big[1_{A_n(t)}\,\Big(|X_n(t,x)|^p+|X_n^-(t,x)|^p\Big)\Big]<\infty\,, \tag{2.29}$$

*and*

$$\sup_{(t,x)\in[0,T]\times\mathbb{R}^2}\Big\|\Big(X_n(t,x)-X_n^-(t,x)\Big)\,1_{A_n(t)}\Big\|_p \le C\,n\,2^{-n\,\frac{1+\beta}{2}}\,. \tag{2.30}$$

*Proof*. We want to show that, for any $p\in[1,\infty)$,

$$\sup_{n\ge1}\sup_{k\ge0}\sup_{(t,x)\in[0,T]\times\mathbb{R}^2} E\Big[1_{A_n(t)}\,\Big(|X_n^k(t,x)|^p+|\bar{X}_n^k(t,x)|^p\Big)\Big]<+\infty\,. \tag{2.31}$$

Indeed, (2.29) is a consequence of (2.31) and (2.28). For $r \le t$ consider the decomposition

$$E\Big(|X_n^{k+1}(t,r;x)|^p\ 1_{A_n(t)}\Big) \le C\ \sum_{i=1}^{6}\ T_n^{k+1,i}(t,r;x)\,, \tag{2.32}$$

where

$$\begin{aligned}
T_n^{k+1,1}(t,r;x) &= |X^0(t,x)|^p,\\
T_n^{k+1,2}(t,r;x) &= E\Big(\Big|\int_0^r\int_{\mathbb{R}^2} S(t-s,\,x-y)\ A(X_n^k(s,y))\ F(ds,\,dy)\Big|^p\ 1_{A_n(t)}\Big),\\
T_n^{k+1,3}(t,r;x) &= E\Big(\Big|\langle S(t-\cdot,x-*)\,B(\bar X_n^k(\cdot,*))\ 1_{[0,r]}(\cdot),\ \omega^n\rangle_H\Big|^p\ 1_{A_n(t)}\Big),\\
T_n^{k+1,4}(t,r;x) &= E\Big(\Big|\langle S(t-\cdot,x-*)\ [B(X_n^k)-B(\bar X_n^k)]\,(\cdot,*)\ 1_{[0,r]}(\cdot)\,,\omega^n\rangle_H\Big|^p\ 1_{A_n(t)}\Big),\\
T_n^{k+1,5}(t,r;x) &= E\Big(\Big|\langle S(t-\cdot,x-*)\,D(X_n^k)\,(\cdot,*)\ 1_{[0,r]}(\cdot),\ h\rangle_H\Big|^p\ 1_{A_n(t)}\Big),\\
T_n^{k+1,6}(t,r;x) &= E\Big(\Big|\int_0^r\int_{\mathbb{R}^2} S(t-s,\,x-y)\ b(X_n^k(s,y))\,ds\,dy\Big|^p\ 1_{A_n(t)}\Big).
\end{aligned}$$

Under hypotheses weaker than (C2), we have proved in Millet and Sanz-Solé (1997):

$$|X^0(t,x)| \le C\Big(\|v_0\|_{q_0} + \|\nabla u_0\|_{q_0}\Big). \tag{2.33}$$

Burkholder's and Hölder's inequalities yield

$$T_n^{k+1,2}(t,r;x) \le C\ \int_0^r J(t-s)\ \Big[1+\sup_{y\in\mathbb{R}^2}\ E(|X_n^k(s,y)|^p\ 1_{A_n(s)})\Big]\,ds\,. \tag{2.34}$$

Let $\tau_n$ be the operator defined on real-valued functions as follows:

$$\tau_n(\rho)\ (s,x) = \rho\Big((s+2^{-n})\wedge T,\ x\Big).$$

Let $\mathcal{E}_n$ be the closed subspace of $H$ generated by the orthonormal system

$$2^n\ T^{-1}\ 1_{\Delta_i}(\cdot)\otimes e_j(*),\ \ i=0,\ldots,\ 2^n-1,\ \ j=1,\ldots,\,n\,,$$

and denote by $\pi_n$ the orthogonal projection operator on $\mathcal{E}_n$. Then

$$\begin{aligned}
T_n^{k+1,3}(t,r;x) = E\Big(\Big|\int_0^1\int_{\mathbb{R}^2}(\pi_n\circ\tau_n)\ \big[S(t-\cdot,\,x-*)\ B(\bar X_n^k(\cdot,*))\\
\times\ 1_{[0,r]}(\cdot)\,1_{A_n(\cdot)}\big]\ (s,y)\ F(ds,\,dy)\Big|^p\Big)
\end{aligned}$$

and, by Burkholder's and Hölder's inequalities, if $J$ is defined by (4.1),

$$\begin{aligned}
T_n^{k+1,3}(t,r;x) &\le C\ E\|(\pi_n\circ\tau_n)\ (S(t-\cdot,\,x-*)\ B(\bar X_n^k(\cdot,*)\ 1_{[0,r]}(\cdot)\ 1_{A_n(\cdot)})\|_H^p\\
&\le C\int_0^r J(t-s)\ \Big[1+\sup_{y\in\mathbb{R}^2}\ E(|\bar X_n^k(s,y)|^p\ 1_{A_n(s)})\Big]ds\,. \qquad (2.35)
\end{aligned}$$

Schwarz's and Hölder's inequalities, (2.12) and (2.21) imply

$$\begin{aligned} T_n^{k+1,4}(t,r;x) &\le E\Big( \|\omega^n \, 1_{A_n(t)}\|_H^p \, \| \, S(t-\cdot,\, x-*) \\ &\qquad \times \, [B(X_n^k) - B(\bar{X}_n^k)]\,(\cdot,*)\; 1_{[0,r]}(\cdot)\; 1_{A_n(t)}\|_H^p \Big) \\ &\le C\, n^p\, 2^{\frac{np}{2}} \int_0^r J(t-s) \sup_{y\in\mathbb{R}^2} E(\, |X_n^k(s,y) - \bar{X}_n^k(s,y)|^p)\, ds \\ &\le C\, n^{2p}\, 2^{-n\beta\,\frac{np}{2}} \\ &\qquad \times \int_0^r J(t-s)\Big[1 + \sup_{(u,y)\in[0,s]\times\mathbb{R}^2} E(\, |X_n^{k-1}(u,y)|^p\, 1_{A_n(r)})\Big] ds\,. \end{aligned} \tag{2.36}$$

Analogously,

$$T_n^{k+1,5}(t,r;x) \le C\, \|h\|_H^p \int_0^r J(t-s)\, \Big[1 + \sup_{y\in\mathbb{R}^2} E(\, |X_n^k(s,y)\,|^p\, 1_{A_n(s)})\Big] ds\,, \tag{2.37}$$

$$T_n^{k+1,6}(t,r;x) \le C \int_0^r (t-s)\, \Big[1 + \sup_{y\in\mathbb{R}^2} E(\, |X_n^k(s,y)\,|^p\, 1_{A_n(s)})\Big]\, ds\,. \tag{2.38}$$

Therefore, (2.32)-(2.38) yield

$$\begin{aligned} E\Big(|X_n^{k+1}(t,r;x)|^p\, 1_{A_n(t)}\Big) &\le \int_0^r \Big[1 + \sup_{(u,y)\in[0,s]\times\mathbb{R}^2} \Big\{ E\big(\, \big[|X_n^k(u,y)|^p \\ &\qquad + |X_n^{k-1}(u,y)|^p + |\bar{X}_n^k(u,y)|^p\big]\; 1_{A_n(s)}\big)\Big\}\Big]\, ds\,. \end{aligned} \tag{2.39}$$

Set, for any $k \ge 0\,,\ t \in [0,T]$,

$$\varphi_n^k(t) = \sup_{(s,y)\in[0,t]\times\mathbb{R}^2} E\Big(\big(\, |X_n^k(s,y)|^p + |\bar{X}_n^k(s,y)|^p\big)\; 1_{A_n(s)}\Big).$$

Then, using (2.39) with $r = t$ and $r = t_n$ and adding term by term the corresponding inequalities, we obtain

$$\varphi_n^{k+1}(t) \le C \int_0^t [1 + \varphi_n^k(s) + \varphi_n^{k-1}(s)]\, ds\,, \tag{2.40}$$

with the convention $\varphi_n^{-1}(\cdot) \equiv 0$. Since by (2.33)

$$\varphi_n^0(t) \le 2 \sup_{x\in\mathbb{R}^2} |X^0(t,x)|^p \le C\,,$$

(2.40) yields

$$\sup_{n\ge1}\ \sup_{k\ge0}\ \sup_{t\in[0,T]} \varphi_n^k(t) \le C\,,$$

which establishes (2.31). Finally, the inequalities (2.22) and (2.29) imply (2.30), which completes the proof of the lemma. □

In the sequel $K$ denotes an arbitrary compact subset of $\mathbb{R}^2$. For any $s, t, \bar{t} \in [0,T],\ x, \bar{x} \in K,\ y \in \mathbb{R}^2$, set

$$\gamma(t,\bar{t},x,\bar{x};\, s,y) = S(t-s,\, x-y) - S(\bar{t}-s,\, \bar{x}-y)$$

and

$$\Gamma(t,\bar t,x,\bar x;\ s,y)=|\gamma(t,\bar t,x,\bar x;\ s,y)|\,.$$

**Lemma 2.8.** *Assume that the function $f$ satisfies the condition (C1). For any* $0<\gamma<\frac{\beta}{2(1+\beta)}$ , $t,\bar t\in[0,T],\ x,\bar x\in K$ ,

$$\|\Gamma(t,\bar t,x,\bar x;\ \cdot,*)\|_H\le C\,(\,|t-\bar t|^\gamma+|x-\bar x|^\gamma)\,, \tag{2.41}$$

$$\int_0^T\int_{\mathbb{R}^2}\Gamma(t,\bar t,x,\bar x;\ s,y)\,ds\,dy\le C\,(\,|t-\bar t|^{\frac12}+|x-\bar x|^{\frac12})\,. \tag{2.42}$$

*Proof.* Assume $0\le t\le\bar t\le T$ and set

$$\begin{aligned}
\Gamma_1(t,\bar t,x,\bar x;\ s,y) &= \Big(S(t-s,\ x-y)-S(\bar t-s,\ x-y)\Big)\ 1_{[0,t]}(s)\,,\\
\Gamma_2(t,\bar t,x,\bar x;\ s,y) &= \Big|\Big(S(\bar t-s,\ x-y)-S(\bar t-s,\ \bar x-y)\Big)\ 1_{[0,t]}(s)\Big|\,,\\
\Gamma_3(t,\bar t,x,\bar x;\ s,y) &= S(\bar t-s,\ \bar x-y)\ 1_{[t,\bar t]}(s)\,;
\end{aligned}$$

then

$$\|\Gamma(t,\bar t,x,\bar x;\ \cdot,*)\|_H^2\le C\ \sum_{i=1}^3\ \|\Gamma_i(t,\bar t,x,\bar x;\ \cdot,*)\|_H^2\,.$$

For $i=1,2$, it is easy to check

$$\|\Gamma_1(t,\bar t,x,\bar x;\ \cdot,*)\|_H^2\le\mu_{t,\bar t-t}+\tilde\mu_{t,\bar t-t}+2\Big(\mu_{t,\bar t-t}\ \tilde\mu_{t,\bar t-t}\Big)^{\frac12},$$

$$\|\Gamma_2(t,\bar t,x,\bar x;\ \cdot,*)\|_H^2\le M_{t,\bar x-x}+N_{t,\bar x-x}+2\Big(M_{t,\bar x-x}\ N_{t,\bar x-x}\Big)^{\frac12},$$

where $\mu_{t,\bar t-t},\ \tilde\mu_{t,\bar t-t},\ M_{t,\bar x-x},\ N_{t,\bar x-x}$ are defined in (4.6)–(4.9), respectively.
Finally,

$$\|\Gamma_3(t,\bar t,x,\bar x;\ \cdot,*)\|_H^2=\mu(\bar t-t)\,.$$

Thus the estimates proved in Lemmas A1 and A5 of Millet and Sanz-Solé (1997) (see (4.12) and (4.13)) show (2.41).
In order to prove (2.42), set $\Gamma(t,\bar t,x,\bar x;\ s,y)\le\sum_{i=1}^2\ \bar\Gamma_i(t,\bar t,x,\bar x;\ s,y)$, with

$$\bar\Gamma_1(t,\bar t,x,\bar x;\ s,y)=\ |\,S(t-s,\ x-y)-S(\bar t-s,\ x-y)\,|\,,$$

$$\bar\Gamma_2(t,\bar t,x,\bar x;\ s,y)=\ |\,S(\bar t-s,\ x-y)-S(\bar t-s,\ \bar x-y)\,|\,.$$

Assume $0\le t\le\bar t$. Then

$$\int_0^T ds\ \int_{\mathbb{R}^2}dy\ \bar\Gamma_1(t,\bar t,x,\bar x;\ s,y)\le C\Big(\nu_{t,\bar t-t}+\tilde\nu_{t,\bar t-t}+\nu(\bar t-t)\Big),$$

with $\nu(\bar t-t),\ \nu_{t,\bar t-t},\ \tilde\nu_{t,\bar t-t}$, defined in (4.3)–(4.5). Hence the inequalities (4.3) and (4.10) imply

$$\int_0^T ds\ \int_{\mathbb{R}^2}dy\ \bar\Gamma_1(t,\bar t,x,\bar x;\ s,y)\le C\ (\bar t-t)^{\frac12}. \tag{2.43}$$

Moreover,

$$\int_0^T ds \int_{\mathbb{R}^2} dy\, \bar\Gamma_2(t,\bar t,x,\bar x;\, s,y) \le \bar\Gamma_{2,1}(t,\bar t,x,\bar x) + 2\bar\Gamma_{2,2}(t,\bar t,x,\bar x)\,,$$

with

$$\begin{aligned}
\bar\Gamma_{2,1}(t,\bar t,x,\bar x) &= \int_0^{\bar t-\frac{|x-\bar x|}{2}} ds \int_{\substack{|x-y|<\bar t-s\\ |\bar x-y|<\bar t-s}} dy\, |S(\bar t-s,\, x-y) - S(\bar t-s,\, \bar x-y)|\,,\\
\bar\Gamma_{2,2}(t,\bar t,x,\bar x) &= \int_0^{\bar t} ds \int_{\substack{|\bar x-y|<\bar t-s\\ |x-y|\ge\bar t-s}} dy\, S(\bar t-s,\, \bar x-y)\,.
\end{aligned}$$

Using (A.24) in Millet and Sanz-Solé (1997, Lemma A4), we obtain $\bar\Gamma_{2,1}(t,\bar t,x,\bar x) \le C\,|x-\bar x|^{\frac12}$. Finally, (1.31) in Millet and Sanz-Solé (1997) implies $\bar\Gamma_{2,2}(t,\bar t,x,\bar x) \le C\,|x-\bar x|^{\frac12}$. Thus

$$\int_0^T ds \int_{\mathbb{R}^2} dy\, \bar\Gamma_2(t,\bar t,x,\bar x;\, s,y) \le C\, |x-\bar x|^{\frac12}. \tag{2.44}$$

The inequalities (2.43) and (2.44) show (2.42) and conclude the proof of the Lemma. ☐

In the next Proposition, we show that the sequence of processes $\{X_n(t,x),\, n\ge 1\}$ satisfies the assumption (P1) of Lemma 4.1. It proves estimates similar to those in Millet and Sanz-Solé (1997, Proposition 1.4) which are uniform in $n$.

**Proposition 2.9.** *Assume (C1), (C2), (C3'). For any $p\in[1,\infty)$, $0\le t\le\bar t\le T$, $x,\bar x\in K$ , $\gamma\in ]0,\frac{\beta}{2(1+\beta)}[$ ,*

$$\sup_n \|\,(X_n(t,x)-X_n(\bar t,\bar x))\ 1_{A_n(\bar t)}\|_p \le C\ (\,|t-\bar t|^\gamma + |x-\bar x|^\gamma)\,.$$

*Proof*. Consider the decomposition

$$E\Big(|X_n(t,x)-X_n(\bar t,\bar x)|^p\ 1_{A_n(\bar t)}\Big) \le C\ \sum_{i=1}^6\ R_n^i(t,\bar t;\, x,\bar x)\,,$$

where

$$\begin{aligned}
R_n^1(t,\bar t;\, x,\bar x) &= |X_0(t,x) - X_0(\bar t,\bar x)|^p,\\
R_n^2(t,\bar t;\, x,\bar x) &= E\Big(\Big|\int_0^1\int_{\mathbb{R}^2} \gamma(t,\bar t,\, x,\bar x;\, s,y)\; A(X_n(s,y))\; F(ds,\, dy)\Big|^p\; 1_{A_n(\bar t)}\Big),\\
R_n^3(t,\bar t;\, x,\bar x) &= E\Big(|\langle\gamma(t,\bar t,\, x,\bar x;\, \cdot,*)\; B(X_n^-(\cdot,\, *)),\; \omega^n\rangle|_H^p\; 1_{A_n(\bar t)}\Big),\\
R_n^4(t,\bar t;\, x,\bar x) &= E\Big(|\langle\gamma(t,\bar t,\, x,\bar x;\, \cdot,*)\; (B(X_n) - B(X_n^-))\,(\cdot,\, *),\; \omega^n\rangle|_H^p\; 1_{A_n(\bar t)}\Big),\\
R_n^5(t,\bar t;\, x,\bar x) &= E\Big(|\langle\gamma(t,\bar t,\, x,\bar x;\, \cdot,*)\; D(X_n(\cdot,\, *)),\; h\rangle|_H^p\; 1_{A_n(\bar t)}\Big),\\
R_n^6(t,\bar t;\, x,\bar x) &= E\Big(\Big|\int_0^1\int_{\mathbb{R}^2} \gamma(t,\bar t,\, x,\bar x;\, s,y)\; b(X_n(s,y))\; ds\, dy\Big|^p\; 1_{A_n(\bar t)}\Big)\;.
\end{aligned}$$

In the proof of Proposition 1.4 in Millet and Sanz-Solé (1997) we have checked that for $\delta = \frac{\beta}{2(1+\beta)}$,

$$R_n^1(t,\bar t;\, x,\bar x) \le C\,(\,|t-\bar t|^\delta + |x-\bar x|^\delta)\,. \tag{2.45}$$

Burkholder's and Hölder's inequalities yield

$$R_n^2(t,\bar t;\, x,\bar x) \le C\; \|\Gamma(t,\bar t,\, x,\bar x;\, \cdot,*)\|_H^p \Big[1 + \sup_{(s,x)\in[0,T]\times\mathbb{R}^2} E(\,|X_n(s,x)|^p\; 1_{A_n(s)})\Big]\,. \tag{2.46}$$

Using the operators $\tau_n$ and $\pi_n$ introduced in the proof of Lemma 2.7 and standard arguments, we obtain

$$\begin{aligned}
R_n^3(t,\bar t;\, x,\bar x) &= E\Big(\Big|\int_0^1\int_{\mathbb{R}^2} (\pi_n\circ\tau_n)\,(\Gamma(t,\bar t,\, x,\bar x;\, \cdot,*)\; B(X_n^-(\cdot,\, *))\\
&\qquad \times\; 1_{A_n(\bar t)})\,(s,y)\; F(ds,\, dy)\Big|^p\Big)\\
&\le C\; \|\Gamma(t,\bar t,\, x,\bar x;\, \cdot,*)\|_H^p\; \Big[1 + \sup_{(s,x)\in[0,T]\times\mathbb{R}^2} E(\,|X_n^-(s,x)|^p\; 1_{A_n(s)})\Big]\,. \qquad (2.47)
\end{aligned}$$

Schwarz's and Hölder's inequalities, (2.12) and (2.30) imply

$$\begin{aligned}
R_n^4(t,\bar t;\, x,\bar x) &\le \Big\{E\Big(\|\omega_n\|_H^{2p}\; 1_{A_n(T)}\Big)\, E\Big(\|\gamma(t,\bar t,\, x,\bar x;\, \cdot,*)\\
&\qquad \times\; (B(X_n) - B(X_n^-))\,(\cdot,\, *)\|_H^{2p}\; 1_{A_n(\bar t)}\Big)\Big\}^{\frac12}\\
&\le C\; n^p\; 2^{n\frac{p}{2}}\; \|\Gamma(t,\bar t,\, x,\bar x;\, \cdot,*)\|_H^p\\
&\qquad \times\Big\{\sup_{(s,x)\in[0,T]\times\mathbb{R}^2} E(\,|X_n(s,x) - X_n^-(s,x)|^{2p}\; 1_{A_n(s)})\Big\}^{\frac12}\\
&\le C\; n^{2p}\; 2^{-n\beta\frac{p}{2}}\; \|\Gamma(t,\bar t,\, x,\bar x;\, \cdot,*)\|_H^p\,. \qquad (2.48)
\end{aligned}$$

Finally,

$$R_n^5(t,\bar t;\,x,\bar x) \le C\,\|h\|_H^p\ \|\Gamma(t,\bar t,\,x,\bar x;\,\cdot,*)\|_H^p \Big[1+\sup_{(s,x)\in[0,T]\times\mathbb{R}^2} E(\,|X_n(s,x)|^p\,1_{A_n(s)})\Big]\,, \tag{2.49}$$

$$\begin{aligned} R_n^6(t,\bar t;\,x,\bar x) \ &\le\ C\,\Big(\int_0^1\int_{\mathbb{R}^2}\Gamma(t,\bar t,\,x,\bar x;\,s,y)\,ds\,dy\Big)^p \\ &\qquad\times\Big[1+\sup_{(s,x)\in[0,T]\times\mathbb{R}^2} E(\,|X_n(s,x)|^p\ 1_{A_n(s)})\Big]. \end{aligned} \tag{2.50}$$

Hence, (2.45)-(2.50), (2.29) and Lemma 2.8 yield the Proposition. □

**Remark 2.10.** Proposition 2.9 establishes the $\gamma$-Hölder continuity for the trajectories of $X_n$ on $A_n(T)$, because the sets $A_n(t),\ t\in[0,T]$, are decreasing. Fix $n\ge 1$. In the proof of Lemma 2.3 we have shown

$$P\big(A_n(T)^c\big)\le C\ \sqrt{n}\ \exp\Big(-n\Big(\frac{\alpha^2}{2}-\ln\ 2\Big)\Big)\,,$$

for any $\alpha>(2\ ln\,2)^{\frac12}$. Consequently, $\lim_{\alpha\to\infty}\ P(A_n(T))=1$, so that the trajectories of $X_n$ are a.s. $\gamma$-Hölder continuous for any $\gamma<\frac{\beta}{2(1+\beta)}$.

We now prove that the processes $\{X_n(t,x),\ n\ge 1\}$ satisfies the condition (P2) of Lemma 4.1.

**Proposition 2.11.** *Suppose that the conditions (C1), (C2) and (C3') are satisfied. Then, for any $p\in[1,\infty),\ (t,x)\in[0,T]\times K$,*

$$\lim_{n\to\infty}\ \|\,(X_n(t,x)-X(t,x))\ 1_{A_n(t)}\|_p=0\,. \tag{2.51}$$

*Proof.* Set

$$X_n(t,x)-X(t,x)=\sum_{i=1}^8\ U_n^i(t,x)$$

with

$$\begin{aligned}
U_n^1(t,x) &= \int_0^t \int_{\mathbb{R}^2} S(t-s,\, x-y)\ [\,(A+B)\ (X_n(s,y)) \\
&\quad - (A+B)\ (X(s,y))]\ F(ds,\, dy)\,, \\
U_n^2(t,x) &= \langle S(t-\cdot,\, x-*)\ [D(X_n(\cdot,*)) - D(X(\cdot,*))],\ h\rangle_H, \\
U_n^3(t,x) &= \int_0^t \int_{\mathbb{R}^2} S(t-s,\, x-y)\ [\,b(X_n(s,y)) - b(X(s,y))]\, ds\, dy\,, \\
U_n^4(t,x) &= \langle S(t-\cdot,\, x-*)\ [B(X_n(\cdot,*)) - B(X_n^-(\cdot,*))],\ \omega^n\rangle_H, \\
U_n^5(t,x) &= \langle S(t-\cdot,\, x-*)\ [B(X_n^-(\cdot,*)) - B(X^-(\cdot,*))],\ \omega^n\rangle_H, \\
U_n^6(t,x) &= \langle S(t-\cdot,\, x-*)\ B(X^-(\cdot,*)),\ \omega^n\rangle_H \\
&\quad - \int_0^t \int_{\mathbb{R}^2} S(t-s,\, x-y)\ B(X^-(s,y))\ F(ds,\, dy)\,, \\
U_n^7(t,x) &= \int_0^t \int_{\mathbb{R}^2} S(t-s,\, x-y)\ [\,B(X^-(s,y)) - B(X_n^-(s,y))]\, F(ds\, dy)\,, \\
U_n^8(t,x) &= \int_0^t \int_{\mathbb{R}^2} S(t-s,\, x-y)\ [\,B(X_n^-(s,y)) - B(X_n(s,y))]\, F(ds\, dy)\,,
\end{aligned}$$

with $X_n^-$ and $X^-$ defined in (2.16) and (2.17), respectively.
Burkholder's and Hölder's inequalities imply

$$\|U_n^1(t,x)\ 1_{A_n(t)}\|_p^p \le C\ \int_0^t J(t-s)\ \sup_{x\in\mathbb{R}^2}\ \|(X_n(s,x) - X(s,x))\ 1_{A_n(s)}\|_p^p\, ds\,.$$

Schwarz's and Hölder's inequalities yield

$$\|U_n^2(t,x)\ 1_{A_n(t)}\|_p^p \le C\ \|h\|_H^p\ \int_0^t J(t-s)\ \sup_{x\in\mathbb{R}^2}\ \|(X_n(s,x) - X(s,x))\ 1_{A_n(s)}\|_p^p\, ds\,.$$

Analogously,

$$\|U_n^3(t,x)\ 1_{A_n(t)}\|_p^p \le C\ \int_0^t (t-s)\ \sup_{x\in\mathbb{R}^2}\ \|(X_n(s,x) - X(s,x))\ 1_{A_n(s)}\|_p^p\, ds\,.$$

Since

$$U_n^5(t,x) = \int_0^t \int_{\mathbb{R}^2} (\pi_n \circ \tau_n) \Big[S(t-\cdot,\, x-*)\ [B(X_n^-) - B(X^-)]\,(\cdot,*)\ 1_{A_n(\cdot)}\Big](s,y)\ F(ds,\, dy)\,,$$

Burkholder's and Hölder's inequalities easily yield

$$\|\,(U_n^5(t,x) + U_n^7(t,x))\ 1_{A_n(t)}\|_p^p \le C\ \int_0^t J(t-s)\ \sup_{x\in\mathbb{R}^2}\ \|(X_n^-(s,x) - X^-(s,x))\ 1_{A_n(s)}\|_p^p\, ds\,.$$

Thus, (2.30) and Lemma 2.5 ensure

$$\begin{aligned}
\|\,(U_n^5(t,x) + U_n^7(t,x))\ 1_{A_n(t)}\|_p^p \ &\le\ C\ n^p\ 2^{-n(1+\beta)\,\frac{p}{2}} \\
&\quad + C\ \int_0^t J(t-s)\ \sup_{x\in\mathbb{R}^2} \|(X_n(s,x) - X(s,x))\, 1_{A_n(s)}\|^p\ ds\,.
\end{aligned}$$

Thus, by Gronwall's lemma, if suffices to check

$$\sup_{0\le t\le T}\ \sup_{x\in K}\ \|U_n^i(t,x)\ 1_{A_n(t)}\|_p \xrightarrow[n\to\infty]{} 0,\ i=4,6,8\,. \tag{2.52}$$

Schwarz's and Hölder's inequalities, (2.12) and (2.30) imply

$$\begin{aligned}\|U_n^4(t,x)\ 1_{A_n(t)}\|_p^p &\le C\Big\{E\Big(\ \|\omega^n\ 1_{A_n(t)}\|_H^{2p}\Big)\Big\}^{\frac12}\\ &\quad\times \sup_{(t,x)\in[0,T]\times\mathbb{R}^2}\ \Big\{E\Big(\ |X_n(t,x)-X_n^-(t,x)|^{2p}\ 1_{A_n(t)}\Big)\Big\}^{\frac12}\\ &\le C\ n^{2p}\ 2^{-n\beta\frac{p}{2}}.\end{aligned}$$

Burkholder's, Hölder's inequalities and (2.30) imply

$$\begin{aligned}\|U_n^8(t,x)\ 1_{A_n(t)}\|_p^p &\le C\int_0^t J(t-s)\ \sup_{x\in\mathbb{R}^2}\ \Big(\|(X_n^-(s,x)-X_n(s,x))\ 1_{A_n(s)}\|_p^p\Big)ds\\ &\le C\ n^p\ 2^{-n(1+\beta)\frac{p}{2}}\ ;\end{aligned}$$

thus (2.52) holds for $i=4,8$.
Set $\|U_n^6(t,x)\ 1_{A_n(t)}\|_p^p\le C\Big(\bar U_n^{6,1}(t,x)+\bar U_n^{6,2}(t,x)\Big)$, with

$$\begin{aligned}\bar U_n^{6,1}(t,x) &= \Big\|\int_0^t\int_{\mathbb{R}^2}\Big[(\pi_n\circ\tau_n)\Big(S(t-\cdot,\ x-*)\ B(X^-(\cdot,*))\ 1_{A_n(\cdot)}\Big)\\ &\qquad -\pi_n\Big(S(t-\cdot,\ x-*)\ B(X^-(\cdot,*))\ 1_{A_n(\cdot)}\Big)\Big](s,y)\ F(ds,\,dy)\Big\|_p^p,\\ \bar U_n^{6,2}(t,x) &= \Big\|\int_0^t\int_{\mathbb{R}^2}\Big[\pi_n\Big(S(t-\cdot,\ x-*)\ B(X^-(\cdot,*))\ 1_{A_n(\cdot)}\Big)(s,y)\\ &\qquad -S(t-s,\ x-y)\ B(X^-(s,y))\ 1_{A_n(s)}\Big]\ F(ds,\,dy)\Big\|_p^p.\end{aligned}$$

By Burkholder's inequality,

$$\begin{aligned}\bar U_n^{6,1}(t,x) &\le C\ E\Bigg(\Big\|\tau_n\Big(S(t-\cdot,\ x-*)\ B(X^-(\cdot,*))\ 1_{A_n(\cdot)}\Big)\\ &\qquad -S(t-\cdot,\ x-*)\ B(X^-(\cdot,*))\ 1_{A_n(\cdot)}\Big\|_H^p\Bigg).\end{aligned}$$

Define

$$\begin{aligned}\bar U_n^{6,1,1}(t,x) &= E\Big(\Big\|\Big[\tau_n\Big(S(t-\cdot,\ x-*)\Big)-\ S(t-\cdot,\ x-*)\Big]\\ &\qquad\times\ \tau_n\Big(B(X^-(\cdot,*))\ 1_{A_n(\cdot)}\Big)\Big\|_H^p\Big),\\ \bar U_n^{6,1,2}(t,x) &= E\Big(\Big\|S(t-\cdot,\ x-*)\Big[\tau_n\Big(B(X^-(\cdot,*))\ 1_{A_n(\cdot)}\Big)\\ &\qquad -\ B(X^-(\cdot,*))\ 1_{A_n(\cdot)}\Big]\Big\|_H^p\Big).\end{aligned}$$

The process $\{X^-(t,x),\ (t,x) \in [0,T] \times \mathbb{R}^2\}$ defined by (2.17) depends on $n$, but Lemma 2.5 and (2.7) imply

$$\sup_{n\geq 1} \sup_{(t,x)\in[0,T]\times\mathbb{R}^2} E\left(|X^-(t,x)|^p\right) < \infty,\ p \in [1,\infty)\,.$$

This property together with (4.12) yield

$$\bar{U}_n^{6,1,1}(t,x) \leq C\ 2^{-n\gamma p} \tag{2.53}$$

for any $\gamma < \frac{\beta}{2(1+\beta)}$.

Moreover, since $X$ is a particular case of $X_n$, by Lemma 2.5 and Proposition 2.9 we have

$$\sup_{(s,x)\in[0,T]\times\mathbb{R}^2} \|(X^-(s+2^{-n},x) - X^-(s,x))\ 1_{A_n(s)}\|_p \leq C\ 2^{-n\gamma}$$

for $0 < \gamma < \frac{\beta}{2(1+\beta)}$. Consequently,

$$\bar{U}_n^{6,1,2}(t,x) \leq C\ 2^{-n\gamma p}\,. \tag{2.54}$$

The inequalities (2.53) and (2.54) imply $\bar{U}_n^{6,1}(t,x) \leq C\ 2^{-n\gamma p}$ and therefore

$$\lim_{n\to\infty} \sup_{(t,x)\in[0,T]\times K} \bar{U}_n^{6,1}(t,x) = 0\,. \tag{2.55}$$

Let $I_H$ denote the identity operator on the Hilbert space $H$. Burkholder's inequality yield $\bar{U}_n^{6,2}(t,x) \leq C\ Z_n(t,x)$, with, for $(t,x) \in [0,T] \times K$

$$Z_n(t,x) = E\Big(\ \|(\pi_n - I_H)\ (\ S(t-\cdot,\ x-*)\ B(X^-(\cdot,*))\ 1_{A_n(\cdot)})\|_H^p\Big).$$

The sequence $\Big(\ \|(\pi_n - I_H)\ (\ S(t-\cdot,\ x-*)\ B(X^-(\cdot,*)) \times 1_{A_n(\cdot)})\|_H\Big)$, $n \geq 1$ decreases to 0 as $n \to \infty$. It is bounded by $\sup_n\ 2\,\|S(t-\cdot,\ x-*)\ B(X^-(\cdot,*))\|_H$ ; we prove that this last random variable belongs to $L^p(\Omega)$. Indeed, Schwarz's inequality implies

$$E\Big(\sup_n\ \|S(t-\cdot,\ x-*)\ B(X^-(\cdot,*))\|_H^p\Big) \leq C\ \sum_{i=1}^4\ T_i\,,$$

with

$$\begin{aligned}
T_1 &= E\Big(\|S(t-\cdot,\ x-*)\ (1+|X^0(\cdot,*))|\ \|_H^p\Big),\\
T_2 &= E\Big(\sup_n\ \Big|\int_0^t ds \iint dy\,dz\ S(t-s,x-y)\ f(|y-z|)\ S(t-s,x-z)\\
&\qquad \times \Big|\int_0^{s_n}\int_{\mathbb{R}^2} S(s-r,y-\eta)\ (A+B)\ (X(r,\eta))\ F(dr\,,d\eta)\Big|\\
&\qquad\quad \times \Big|\int_0^{s_n}\int_{\mathbb{R}^2} S(s-r,z-\zeta)\ (A+B)\ (X(r,\zeta))\ F(dr\,,d\zeta)\Big|\Big|^{\frac{p}{2}}\Big),
\end{aligned}$$

$$\begin{aligned} T_3 &= E\Big(\sup_n \Big|\int_0^t ds \iint dy\, dz\; S(t-s,x-y)\; f(|y-z|)\; S(t-s,x-z) \\ &\qquad \times |\langle 1_{[0,s_n]}(\cdot)\, S(s-\cdot, y-*)\; D(X(\cdot,*)), h\rangle_H| \\ &\qquad \times |\langle 1_{[0,s_n]}(\cdot)\, S(s-\cdot, z-*)\; D(X(\cdot,*)), h\rangle_H|\Big|^{\frac{p}{2}}\Big), \\ T_4 &= E\Big(\sup_n \Big|\int_0^t ds \iint dy\, dz\; S(t-s,x-y)\; f(|y-z|)\; S(t-s,x-z) \\ &\qquad \times |\int_0^{s_n}\int_{\mathbb{R}^2} S(s-r,y-\eta)\; b(X(r,\eta))\, dr\, d\eta| \\ &\qquad \times |\int_0^{s_n}\int_{\mathbb{R}^2} S(s-r,z-\zeta)\; b(X(r,\zeta))\, dr\, d\zeta|\Big|^{\frac{p}{2}}\Big). \end{aligned}$$

Hölder's inequality implies that

$$\begin{aligned} T_2 &\le C\int_0^t ds \iint dy\, dz\; S(t-s,x-y)\; f(|y-z|)\; S(t-s,x-z) \\ &\qquad \times E\Big(\sup_n \Big|\int_0^{s_n}\int_{\mathbb{R}^2} S(s-r,y-\eta)\; (A+B)\,(X(r,\eta))\, F(dr\,,d\eta)\Big|^p\Big)^{\frac{1}{2}} \\ &\qquad \times E\Big(\sup_n \Big|\int_0^{s_n}\int_{\mathbb{R}^2} S(s-r,z-\zeta)\; (A+B)(X(r,\zeta))\, F(dr\,,d\zeta)\Big|^p\Big)^{\frac{1}{2}}. \end{aligned}$$

Doob's maximal inequality applied to the martingale

$$\Big(\int_0^\tau\int_{\mathbb{R}^2} S(s-r,\, y-\eta)\,(A+B)\,(X(r,\eta))\; F(dr\, d\eta)\,,\; \mathcal{F}_\tau\Big)$$

implies that

$$T_2 \le C\Big[1+\sup_{(s,y)\in[0,T]\times\mathbb{R}^2} E(|X(s,y)|^p)\Big].$$

A similar easier computation using Schwarz's and Hölder's inequality yields

$$T_3+T_4 \le C\;\Big[1+\sup_{(s,y)\in[0,T]\times\mathbb{R}^2} E(|X(s,y)|^p)\Big].$$

Therefore, these estimations and (2.7) imply for $p\in[1,\infty[$,

$$E\Big(\sup_n\; \|S(t-\cdot,\, x-*)\; B(X^-(\cdot,*))\|_H^p\Big) < \infty\,.$$

Thus, by dominated convergence, the sequence $(Z_n(t,x))_{n\ge1}$ decreases to 0. Moreover, $Z_n(t,x)$ is jointly continuous in $(t,x)$; consequently by Dini's Theorem

$$\sup_{(t,x)\in[0,T]\times K} Z_n(t,x)\;\downarrow\; 0 \quad\text{as}\quad n\to\infty\,.$$

Thus,

$$\sup_{(t,x)\in[0,T]\times K} \bar U_n^{6,2}(t,x) \underset{n\to\infty}{\longrightarrow} 0\,. \tag{2.56}$$

The convergences (2.55) and (2.56) complete the proof of (2.52) for $i = 6$, and hence that of the Proposition. □

*Proof of Proposition 2.2.* Since equation (2.2) defining the process $\{X(t,x),\ (t,x) \in [0,T] \times \mathbb{R}^2\}$ is a particular case of equation (2.1) which defines $\{X_n(x,t),\ (t,x) \in [0,T] \times \mathbb{R}^2\}$ , Propositions 2.9 and 2.11 ensure the validity of conditions (P1) and (P2) of Lemma 4.1 for the sequence of processes $Y_n(t,x) := X_n(t,x) - X(t,x)$ and the sequence of adapted sets $B_n(t) := A_n(t)$ defined in (2.11). Therefore, given any $0 < \gamma < \frac{\beta}{2(1+\beta)}$ , $p \in [1,\infty)$,

$$\lim_{n\to\infty} E\Big( \|X_n - X\|_{\gamma,K}^p \ 1_{A_n(T)} \Big) = 0\,, \tag{2.57}$$

where $\|\cdot\|_{\gamma,K}$ is given by (1.5).
Fix $\varepsilon > 0$; by Lemma 2.3 there exists $n_0 \in \mathbb{N}$ be such that $P(A_n(T)^c) < \varepsilon$. Then, for any $\eta > 0$,

$$\begin{aligned} P\big(\|X_n - X\|_{\gamma,K} > \eta\big) &\le \varepsilon + P\big(\|X_n - X\|_{\gamma,K} > \eta\,,\ A_n(T)\big) \\ &\le \varepsilon + \eta^{-p}\, E(\, \|X_n - X\|_{\gamma,K}^p \ 1_{A_n(T)}). \end{aligned} \tag{2.58}$$

Since $\varepsilon > 0$ is arbitrary, (2.57) and (2.58) show (2.4). □

## 3. Approximation in $L^p$

In the previous section, we have proved an approximation theorem in probability, by showing the $L^p$ convergence of the sequence $X_n$ localized by $A_{n,\,M(n)}$. The aim of this section is to check that under a stronger growth assumption on the coefficients, a slight modification of the proof yields the $L^p$-convergence of $X^n$ to $X$ without localization. Let us introduce the following growth condition:

(C4') There exists $\delta \in (0,1)$ and a constant $C > 0$ such that for $x \in \mathbb{R}^2$,

$$|A(x)| + |B(x)| + |D(x)| + |b(x)| \le C\ (1 + |x|^\delta)\,.$$

Then we have the following

**Proposition 3.1.** *Assume (C1), (C2) and (C4'), and let $X$ and $X^n$ be defined by (2.1) and (2.2) respectively. For any $\gamma \in \Big(0,\ \frac{\beta}{2(1+\beta)}\Big)$, every compact subset $K \subset \mathbb{R}^2$ and every $p \in [1,+\infty)$,*

$$\lim_n \Big\| \|X_n - X\|_{\gamma,K} \Big\|_p = 0 \tag{3.1}$$

The proof is very similar to that of Proposition 2.2, and will only be sketched. It depends on several technical lemmas, which are "unlocalized" versions of Lemmas 2.6 and 2.7.

**Lemma 3.2.** *Suppose that the conditions (C1), (C2) and (C4') hold. Then for any* $p \in [1, +\infty)$ , $t \in [0, T]$ ,$\delta' \in ]\delta, 1[$ *and* $n \geq 1$,

$$\sup_{(s,y)\in[0,t]\times\mathbb{R}^2} E(\,|X_n^k(s,y) - \bar{X}_n^k(s,y)|^p)$$

$$\leq C\; n^{\frac{p}{2}}\; 2^{-n\,\frac{1+\beta}{2}\,p}\,\Big[1 + \sup_{(s,y)\in[0,t]\times\mathbb{R}^2} E(\,|X_n^{k-1}(s,y)|^{\delta' p})\Big]\,, \tag{3.2}$$

*and*

$$\sup_{(s,y)\in[0,t]\times\mathbb{R}^2} E(\,|X_n(s,y) - X_n^-(s,y)|^p)$$

$$\leq C\; n^{\frac{p}{2}}\; 2^{-n\,\frac{1+\beta}{2}\,p}\,\Big[1 + \sup_{(s,y)\in[0,t]\times\mathbb{R}^2} E(\,|X_n(s,y)|^{\delta p})\Big]. \tag{3.3}$$

*Proof*. Consider the decomposition

$$E(\,|X_n^k(t,x) - \bar{X}_n^k(t,x)|^p) \leq C\,\sum_{i=1}^4 \tilde{T}^{k,i}\,,$$

where each term $\tilde{T}_n^{k,i}$ is deduced from the corresponding term $T_n^{k,i}$ introduced in (2.23) by removing $1_{A_n(t)}$.
Let $\bar{p}$ and $\bar{q}$ be conjugate exponents such that $\delta\bar{p} = \delta'$. Then Schwarz's and Hölder's inequalities, (1.14) and (4.11) imply

$$\begin{aligned}
\tilde{T}^{k,2}(t,x) \;&\leq\; E\Big(\,\|1_{(t_n,t]}\;\omega_n\|_H^{\bar{q}p}\Big)^{\frac{1}{\bar{q}}}\\
&\qquad\times\,E\Big(\|1_{(t_n,t]}(\cdot)\;S(t-\cdot,\,x-*)\;B(X_n^{k-1}(\cdot,*))\|_H^{\bar{p}p}\Big)^{\frac{1}{\bar{p}}}\\
&\leq\; C\;n^{\frac{p}{2}}\;2^{n\,\frac{p}{2}}\;2^{-n\,\frac{p}{2}}\\
&\qquad\times\,\Big[\mu(t-t_n)\,\Big(1 + \sup_{(s,y)\in[0,t]\times\mathbb{R}^2} E(\,|X_n^{k-1}(s,y)|^{\delta\bar{p}p})\Big)\Big]^{\frac{1}{\bar{p}}}\\
&\leq\; C\;n^{\frac{p}{2}}\;2^{-n\,\frac{(1+\beta)}{2}\,p}\,\Big[1 + \sup_{(s,y)\in[0,t]\times\mathbb{R}^2} E(\,|X_n^{k-1}(s,y)|^{\delta' p})\Big].
\end{aligned}$$

The upper estimates of $\tilde{T}_n^{k,i}$, $i = 1, 3, 4$ are obtained by means of a straightforward modification of that of $T_n^{k,i}$ in the proof of Lemma 2.6; this concludes the proof of (3.2).
Using the arguments in the proof of Millet and Sanz-Solé (1997, Theorem 1.2), we obtain the convergence of the Picard iteration scheme, i.e., for $p \in [1, +\infty)$,

$$\lim_{k\to\infty}\;\sup_{(s,x)\in[0,T]\times\mathbb{R}^2}\Big(\|X_n^k(s,x) - X_n(s,x)\|_p + \|\bar{X}_n^k(s,x) - X_n^-(s,x)\|_p\Big) = 0\,. \tag{3.4}$$

Therefore, (3.2) and (3.4) yield (3.3). □

We now prove $L^p$ convergence of $X_n^-(s,y)$ to $X_n(s,y)$.

**Lemma 3.3.** *Assume (C1), (C2) and (C4'); then for* $p \in [1, +\infty)$*,*

$$\sup_{n\geq 1} \sup_{(t,x)\in[0,T]\times\mathbb{R}^2} \Big( \|X_n(t,x)\|_p + \|X_n^-(t,x)\|_p \Big) < \infty \,, \tag{3.5}$$

*and*

$$\sup_{(t,x)\in[0,T]\times\mathbb{R}^2} \|X_n(t,x) - X_n^-(t,x)\|_p \leq C\; n^{\frac12}\; 2^{-n\,\frac{1+\beta}{2}} \,. \tag{3.6}$$

*Proof*. The proof reduces to that of

$$\sup_{n\geq 1} \sup_{k\geq 0} \sup_{(t,x)\in[0,T]\times\mathbb{R}^2} \Big( \|X_n^k(t,x)\|_p + \|\bar X_n^k(t,x)\|_p \Big) < \infty \,. \tag{3.7}$$

Indeed, (3.4) and (3.7) imply (3.5), while (3.3) and (3.5) yield (3.6). For $r \leq t$, consider the decomposition

$$E\Big( |X_n^{k+1}(t,r;\, x)|^p \Big) \leq C \sum_{i=1}^{6} \tilde T_n^{k+1,i}\; (t,r;\, x) \,, \tag{3.8}$$

where $\tilde T_n^{k+1,i}\;(t,r;\,x)$ is deduced from the term $T_n^{k+1,i}\;(t,r;\,x)$ in (2.32) by removing $1_{A_n(t)}$. The arguments used to upper estimate $\tilde T_n^{k+1,i}$ for $i \neq 4$ are similar to that in Lemma 2.7 and are omitted. They yield the analogues of (2.33)–(2.35), (2.37) and (2.38).
Let $\delta' \in ]\delta, 1[$ , $\bar p = (\delta')^{-1}$, $\bar p$ and $\bar q$ be conjugate exponents; then Schwarz's and Hölder's inequalities and (1.13) yield

$$\begin{aligned}
\tilde T_n^{k+1,4}\;(t,r;\,x) \;&\leq\; E\Big( \|\omega^n\|_H^{\bar q p} \Big)^{\frac{1}{\bar q}} \\
&\qquad \times E\Big( \|S(t-\cdot,\, x-*)\; 1_{[0,r]}(\cdot)\; [B(X_n^k) - B(\bar X_n^k)]\; (\cdot,*)\|_H^{\bar p p} \Big)^{\frac{1}{\bar p}} \\
&\leq\; C\; n^{\frac{p}{2}}\; 2^{n\,\frac{p}{2}} \Big[ \sup_{(s,y)\in[0,r]\times\mathbb{R}^2} E\Big( |(X_n^k - \bar X_n^k)(s,y)|^{\bar p p} \Big) \Big]^{\frac{1}{\bar p}}.
\end{aligned}$$

Hence (3.2) implies

$$\begin{aligned}
\tilde T_n^{k+1,4}\;(t,r;\,x) \;&\leq\; C\; n^{\frac{p}{2}}\; 2^{\frac{np}{2}} \Big\{ n^{\frac{\bar p p}{2}}\; 2^{-n(1+\beta)\,\frac{\bar p p}{2}} \\
&\qquad \times \Big[ 1 + \sup_{(s,y)\in[0,r]\times\mathbb{R}^2} E\Big( |X_n^{k-1}(s,y)|^{\delta' \bar p p} \Big) \Big] \Big\}^{\frac{1}{\bar p}} \\
&\leq\; C\; n^p\; 2^{-n\,\beta\,\frac{p}{2}} \Big[ 1 + \sup_{(s,y)\in[0,r]\times\mathbb{R}^2} E\Big( |X_n^{k-1}(s,y)|^p \Big) \Big] \,. \qquad (3.9)
\end{aligned}$$

Set $\tilde\varphi_n^{-1} \equiv 0$, and for every $k \geq 0$,

$$\tilde\varphi_n^k(t) = \sup_{0\leq s\leq t} \sup_{y\in\mathbb{R}^2} E(\, |X_n^k(s,y)|^p + |\bar X_n^k(s,y)|^p) \,.$$

Then for every $k \geq 0$ ,

$$\tilde\varphi_n^{k+1}(t) \leq C \int_0^t \Big[ 1 + \tilde\varphi_n^k(s) + \tilde\varphi_n^{k-1}(s) \Big] ds \,.$$

Since $\sup_{0\le t\le T}\ \tilde{\varphi}_n^0(t) = C < \infty$, this implies (3.7). □

Replacing (2.29) and (2.30) by (3.5) and (3.6) respectively, the arguments in the proofs of Propositions 2.9 and 2.11 yield the following

**Proposition 3.4.** *Assume (C1), (C2) and (C4'). Let $K$ be a compact subset of $\mathbb{R}^2$ and $p \in [1,\infty)$; then,*
*(i) For $0 < \gamma < \frac{\beta}{2(1+\beta)}$ , $0 \le t \le \bar{t} \le T$ , $x, \bar{x} \in K$ ,*

$$\sup_n \|X_n(t,x) - X_n(\bar{t},\bar{x})\|_p + \|X(t,x) - X(\bar{t},\bar{x})\|_p \le C(\,|t-\bar{t}|^\gamma + |x-\bar{x}|^\gamma)\,. \tag{3.10}$$

*(ii) For $(t,x) \in K$,*

$$\lim_n \|X_n(t,x) - X(t,x)\|_p = 0\,. \tag{3.11}$$

*Proof of Proposition 3.1.* To conclude the proof of this proposition, it suffices to apply Lemma A.1 in Bally, Millet and Sanz-Solé (1995). Indeed, the results proved in the previous Proposition ensure the validity of the hypothesis of that Lemma. □

## 4. Appendix

This section quotes some notations introduced in our previous paper (Millet and Sanz-Solé, 1997), which are extensively used along the paper. It also contains a technical result.

For any $t \in [0,T],\ h \geq 0,\ \xi \in \mathbb{R}^2$, set

$$(4.1)\quad J(t) = \int_{|y|<|x|<t} \frac{1}{\sqrt{t^2-|x|^2}}\, f(\,|x-y|\,)\, \frac{1}{\sqrt{t^2-|y|^2}}\, dx\, dy\,,$$

$$(4.2)\quad \mu(t) = \int_0^t ds \int_{\mathbb{R}^2} dx \int_{\mathbb{R}^2} dy\ S(s,x)\ f(\,|x-y|\,)\ S(s,y) = \frac{1}{2\,\pi^2} \int_0^t J(s)\, ds\,,$$

$$(4.3)\quad \nu(t) = \frac{1}{2\,\pi} \int_0^t ds \int_{|x|<s} \frac{dx}{\sqrt{s^2-|x|^2}} = \frac{t^2}{2}\,,$$

$$(4.4)\quad \nu_{t,h} = \int_0^t ds \int_{|y|<s} dy\ \Big(S(s,y) - S(s+h,y)\Big),$$

$$(4.5)\quad \tilde{\nu}_{t,h} = \int_0^t ds \int_{s\leq |y|<s+h} dy\ S(s+h,y)\,,$$

$$\mu_{t,h} = \int_0^t ds \int_{|y|<s} dy \int_{|z|<s} dz\ [S(s,y) - S(s+h,y)]\ f(\,|y-z|\,)$$

$$(4.6)\quad \times\ [S(s,z) - S(s+h,z)]\,,$$

$$(4.7)\quad \tilde{\mu}_{t,h} = \int_0^t ds \int_{s\leq|y|<s+h} dy \int_{s\leq|z|<s+h} dz\ S(s+h,y)\ f(\,|y-z|\,)\ S(s+h,z)\,,$$

$$(4.8)\quad M_{t,\xi} = \int_0^t \int_{\substack{|y|<s\\ |y-\xi|\geq s}} dy \int_{\substack{|z|<s\\ |z-\xi|\geq s}} dz\ S(s,y)\ f(\,|y-z|\,)\ S(s,z)\,,$$

$$N_{t,\xi} = \int_{\frac{|\xi|}{2}}^t ds \int_{\substack{|y|<s\\ |y-\xi|<s}} dy \int_{\substack{|z|<s\\ |z-\xi|<s}} dz\ |S(s,y) - S(s,y-\xi)|\ f(|y-z|)$$

$$(4.9)\quad \times\ |S(s,z) - S(s,z-\xi)|\,.$$

A direct computation shows

$$(4.10)\quad \nu_{t,h} + \tilde{\nu}_{t,h} \leq C\ h^{\frac{1}{2}}.$$

Assume that $f$ satisfies the assumption (C1); then Lemma A.1 in Millet and Sanz-Solé (1997) implies

$$(4.11)\quad J(t) \leq C\ t^{\beta},\quad \mu(t) \leq C\ t^{\beta+1},\quad t \in [0,T]\,,$$

while for $t \in [0,T],\ \ h \vee |\xi| \leq \frac{1}{2}$ and $0 < \delta < \frac{\beta}{1+\beta}$, Lemma A.5 in Millet and Sanz-Solé (1997) shows that

$$(4.12)\quad \mu_{t,h} + \tilde{\mu}_{t,h} \ \leq\ C\ h^{\delta}\,,$$

$$(4.13)\quad M_{t,\xi} + N_{t,\xi} \ \leq\ C\ |\xi|^{\delta}\,.$$

The following lemma is a localized version of Lemma A.1 in Bally, Millet and Sanz-Solé (1995). For the sake of completeness we give the main arguments of the proof.

**Lemma 4.1.** *Let* $\{Y_n(t,x),\ (t,x) \in K_0\},\ n \geq 1$ *be a sequence of processes indexed by* $K_0 = [0,T] \times K$ , $K$ *being a compact set of* $\mathbb{R}^2$.
*Let* $\{B_n(t),\ t \in [0,T]\} \subset \mathcal{F}$ *be a sequence of adapted sets which, for every* $n$*, decreases in* $t$*. Assume that for every* $p \in (1,\infty)$*:*
*(P1) There exists* $\delta > 0$ *such that, for any* $0 \leq t \leq \bar{t} \leq T$ , $x, \bar{x} \in K$ ,

$$\sup_n E\Big( |Y_n(t,x) - Y_n(\bar{t},\bar{x})|^p\, 1_{B_n(\bar{t})}\Big) \leq C\left(|t-\bar{t}| + |x-\bar{x}|\right)^{3+\delta}.$$

*(P2) For every* $(t,x) \in [0,T] \times K$,

$$\lim_{n\to\infty} E\Big(|Y_n(t,x)|^p\, 1_{B_n(t)}\Big) = 0\,.$$

*Then, for any* $\rho \in \left(0, \frac{\delta}{p}\right)$ *and any* $r \in [1,p)$,

$$\lim_{n\to\infty} E\Big(\|Y_n\|^r_{\rho,K}\, 1_{B_n(T)}\Big) = 0\,.$$

*Proof.* Let $\zeta = 2\,d + \delta'$, $d = 3$, $0 < \delta' < \delta$; set $z = (t,x)$, $\bar{z} = (\bar{t},\bar{x})$. Then, by (P1),

$$\int_{K_0}\int_{K_0} E\left(\frac{|Y_n(z) - Y_n(\bar{z})|^p}{|z-\bar{z}|^\zeta}\, 1_{B_n(\bar{t})}\right) dz\, d\bar{z} \leq C\ B',$$

where

$$B' = \int_{K_0}\int_{K_0} |z-\bar{z}|^{-d+\delta-\delta'}\, dz\, d\bar{z} < +\infty\,.$$

Set

$$Z = \int_{K_0}\int_{K_0} \frac{|Y_n(z) - Y_n(\bar{z})|^p}{|z-z'|^\zeta}\, 1_{B_n(\bar{t})}\, dz\, d\bar{z}\,.$$

Clearly, by Fubini's theorem, $E\,(Z) \leq C\ B'$, so that

$$P(Z > \lambda^p) \leq C\lambda^{-p}\ B'.$$

The Garsia-Rodemich-Rumsey Lemma yields

$$|Y_n(z) - Y_n(\bar{z})|\ 1_{B_n(\bar{t})} \leq \bar{C}\ Z^{\frac{1}{p}}\ |z-\bar{z}|^{\rho_0}\,,$$

with $\rho_0 = \frac{\delta'}{p}$. Since $\{B_n(t),\ t \in [0,T]\}$ decreases in $t$, this yields for any $\rho < \frac{\delta}{p}$,

$$\begin{aligned} P\left(\sup_{z\neq\bar{z}} \frac{|Y_n(z) - Y_n(\bar{z})|}{|z-\bar{z}|^\rho} > \lambda,\ B_n(T)\right) &\leq \lambda^{-p}\, E\left(\sup_{z\neq\bar{z}} \frac{|Y_n(z) - Y_n(\bar{z})|^p}{|z-\bar{z}|^{\rho p}}\, 1_{B_n(T)}\right) \\ &\leq C\ \lambda^{-p}\ E(Z) \leq C\ \lambda^{-p}. \end{aligned}$$

Intersecting with the set $B_n(T)$, we now proceed exactly as in Bally, Millet and Sanz-Solé (1995, Lemma A.1) and show that for any $\varepsilon > 0$, $r \in [1,p)$, there exists $N \in \mathbb{N}$ such that for any $n \geq N$,

$$E\Big(\|Y_n\|^r_{\beta,k}\, 1_{B_n(T)}\Big) \leq \varepsilon^r + C\varepsilon\,.$$

Annie Millet, Laboratoire de Probabilités, Université Paris VI, 4, place Jussieu, 75252 PARIS Cedex 05, France

*E-mail address*: amil@ccr.jussieu.fr

Marta Sanz-Solé, Facultat de Matemàtiques, Universitat de Barcelona, Gran Via, 585, 08007 Barcelona, Spain

*E-mail address*: sanz@cerber.mat.ub.es